\documentclass[11pt]{article}
\usepackage{graphicx}
\usepackage{amsmath}
\usepackage{amsfonts}
\usepackage{epsfig}
\usepackage{setspace}

\newcommand{\be}{\begin{equation}}
\newcommand{\ee}{\end{equation}}
\newcommand{\beqn}{\begin{eqnarray}}
\newcommand{\eeqn}{\end{eqnarray}}
\newcommand{\beqns}{\begin{eqnarray*}}
\newcommand{\eeqns}{\end{eqnarray*}}

\newcommand{\Cov}{\mbox{Cov}\ }

\newcommand{\EE}{\ensuremath{{\mathbb E}}}
\newcommand{\II}{\ensuremath{{\mathbb I}}}

\newcommand{\fr}[1]{(\ref{#1})}

\newcommand{\om}{\omega}

\newcommand{\Om}{\Omega}

\newcommand{\bom}{\mbox{\mathversion{bold}$\om$}}

\newtheorem{lemma}{Lemma}
\newtheorem{theorem}{Theorem}

\newtheorem{remark}{Remark}

\begin{document}

\title{\Large{\bf Anisotropic Functional Deconvolution for the irregular design with dependent long-memory errors }}

\author{
\large{ Rida Benhaddou}  \footnote{E-mail address: Benhaddo@ohio.edu}
  \\ \\
Department of Mathematics, Ohio University, Athens, OH 45701} 
\date{Draft Version 2}

\doublespacing
\maketitle
\begin{abstract}
Anisotropic functional deconvolution model is investigated in the bivariate case under long-memory errors when the design points $t_i$, $i=1, 2, \cdots, N$, and $x_l$, $l=1, 2, \cdots, M$, are irregular and follow known densities $h_1$, $h_2$, respectively. In particular, we focus on the case when the densities $h_1$ and $h_2$ have singularities, but $1/h_1$ and $1/h_2$ are still integrable on $[0, 1]$.  Under both Gaussian and sub-Gaussian errors, we construct an adaptive wavelet estimator that attains asymptotically near-optimal convergence rates that deteriorate as long-memory strengthens. The convergence rates are completely new and depend on a balance between the smoothness and the spatial homogeneity of the unknown function $f$, the degree of ill-posed-ness of the convolution operator, the long-memory parameter in addition to the degrees of spatial irregularity associated with $h_1$ and $h_2$. Nevertheless, the spatial irregularity affects convergence rates only when $f$ is spatially inhomogeneous in either direction.


{\bf AMS (2000) Subject Classification: 62G05, 62G20, 62G08 }
 \end{abstract} 

\section{Introduction.}

Consider the problem of estimating a periodic bivariate function $f(. , .)\in L^2([0, 1]^2)$ based on the observations

\be \label{eqirr}
Y(t_i, x_l)=\int^1_0f(s, x_l)g(t-s, x_l)ds + \sigma \varepsilon_{il}, i=1, 2, \cdots, N,\ l=1, 2, \cdots, M.
\ee 
Here, $t_i\neq\frac{i}{N}$ and $x_l\neq \frac{l}{M}$, but $0< t_1< t_2< \cdots <t_N<1$ and $0 < x_1< x_2< \cdots< x_M<1$. In addition, $\varepsilon_{il}$ are zero mean Gaussian or sub-Gaussian random variables that are dependent for different $i$, $i=1, 2, \cdots, N$, but independent for different $l$, $l=1, 2, \cdots, M$ and $\sigma$ is a positive variance constant.   The errors $\varepsilon_{il}$ are independent of the design points $t_i$ and $x_l$. Furthermore, the function $g(. , .)$ is the convolution kernel and is supposed to be known. Suppose that the increasing functions $H_1$ and $H_2$ defined on $[0, 1]$ are known and satisfy $H_1(t_i)=\frac{i}{N}$ and $H_2(x_l)=\frac{l}{M}$. Model \fr{eqirr} is referred to as functional deconvolution and it is motivated by experiments in which one needs to recover a bivariate function using observations  of its convolutions along profiles $x=x_l$, $l=1, 2, \cdots, M$. This situation occurs for example in seismic inversion (see Robinson~(1999)).  Let $\pmb{\varepsilon}^l_N$ be zero mean vector with components $\varepsilon_{il}$,  $i=1, 2, \cdots, N$, and let $\Sigma^l_N=\Cov(\pmb{\varepsilon}^l_N)=\EE\left[\pmb{\varepsilon}^l_N(\pmb{\varepsilon}^l_N)^T\right]$ be its covariance matrix. Before we introduce our first assumption, recall that a random variable $X$ is said to be sub-Gaussian if its distribution is dominated by that of a Gaussian random variable. In particular, the tails of a sub-Gaussian distribution decay at least as fast as those of a Gaussian. In addition, a sub-Gaussian random variable $X$ satisfies the moment condition
$$\|X\|_{\psi_2}:=\sup_{p\geq 1}\{{p^{-1/2}}\left(\EE[|X|^p]\right)^{1/p}\}<\infty.$$
Sub-Gaussian random variables include the Gaussian, Bernoulli or any bounded random variable (see {vershynin}~(2011), Sec. 5.2.3), although sub-Gaussianity represents a class of probability distributions rather than a distribution. \\
  Consider the following assumptions in regards of the vectors $\pmb{\varepsilon}^l_N$ and their covariance matrices $\Sigma^{(l)}_N$, $l=1, 2, \cdots, M$.\\
  \noindent
{\bf Assumption A.0.} The vectors $\pmb{\varepsilon}^l_N$ are of the form 
\be \label{epseta}
\pmb{\varepsilon}^l_N=A^l_N\pmb{\eta}^l_N,
\ee
where $\pmb{\eta}^l_N$ are random vectors with zero-mean independent Gaussian or sub-Gaussian components $\eta^l_{i}$ having variance equal to 1, $i=1, 2, \cdots, N$ such that $\|\eta^l_{i}\|_{\psi_2}< K$, $0 < K< \infty$, and $A^l_N$ are some non-random matrices. Notice that under \fr{epseta}, $\Sigma^l_N=\EE\left[\pmb{\varepsilon}^l_N \left(\pmb{\varepsilon}^l_N\right)^T\right]=\EE\left[A^l_N\pmb{\eta}^l_N\left(\pmb{\eta}^l_N\right)^T\left(A^l_N\right)^T\right]=A^l_N\left(A^l_N\right)^T$.\\
\noindent
{\bf Assumption A.1.  }\label{A1}For each covariance matrix $\Sigma^l_N$, there exist constants $c_{1}$ and $c_{2}$ ($0<c_{1} \leq c_{2} < \infty$), independent of $N$, such that 
\be \label{assum1}
c_{1} N^{1-\alpha}\leq \lambda_{\min}\left(\Sigma^l_N\right)\leq \lambda_{\max}\left(\Sigma^l_N\right)\leq c_{2} N^{1-\alpha}, \ \ 0<\alpha \leq 1,
\ee
where $\alpha$ is the long-memory parameter associated with vector $\varepsilon^l_N$, and $\lambda_{\min}\left(\Sigma^l_N\right)$ and $\lambda_{\max}\left(\Sigma^l_N\right)$ are the smallest and the largest eigenvalues of $\Sigma^l_N$, respectively.\\
\noindent
Assumption $A.1$ is valid for example when $\pmb{\varepsilon}^l_N$ are fractional Gaussian noises or fractional ARIMA, (e.g., see Benhaddou et al.~(2014), Section 2).  

Wavelet deconvolution has been the subject of a large array of papers since the pioneering work of Donoho~(1995). Other relevant articles include Abramovich and  Silverman (1998), Walter  and  Shen (1999), Donoho  and  Raimondo (2004), Johnstone et al.~(2004),  among others.  In the case of functional deconvolution model with $f(t, x)\equiv f(t)$, Pensky and sapatinas~(2009, 2010, 2011) pioneered into the formulation and further development of the problem. 
 
Functional deconvolution model of type \fr{eqirr} with $\alpha=1$ and spatially regular design points, corresponds to the $i.i.d.$ case studied in Benhaddou et al.~(2013) where they constructed an adaptive hard-thresholding wavelet estimator, and showed that it is asymptotically near-optimal under the $L^2$-risk over a wide range of Besov balls. Benhaddou~(2017) extends this work to the case of $L^p$-risk, $1\leq p <\infty$. In these articles, errors are assumed to be $i.i.d.$  Gaussian random variables or Gaussian white noise. However, long-memory (LM) or long-range dependence (LRD) is widely encountered phenomenon, including in geophysical exploration (e.g., see Painter et al.~(1995)),  and thus the importance to extend the work of Benhaddou et al.~(2013) to the LM case. This was addressed in the continuous case by Benhaddou and Liu~(2019), where the noise is assumed to be two-parameter fractional Wiener sheet which is technically associated with zero mean long-memory Gaussian sequences.  

LM has been investigated quite considerably in the standard Fourier deconvolution model and the list includes Wang~(1997), Wishart~(2013), Benhaddou et al.~(2014), Kulik et al.~(2015), Benhaddou~(2016) and most recently in Benhaddou and Liu~(2019). LM has been investigated in the Laplace deconvolution framework in Benhaddou~(2018) when the unknown response function is not periodic, and in the density estimation in Comte, Dedecker and Taupin~(2008).

All of the wavelet-related articles mentioned above assume data to be equispaced. But in practice, this may not be the case. The first to address the wavelet regression estimation under irregularly spaced data points was Cai and Brown~(1998). This was followed by many articles, but we will list only a few and they are in chronological order Chesneau~(2007a, 2007b), Pensky~(2013) and Antoniadis et al.~(2014). 

The objective of the paper is to look into the bivariate functional deconvolution in the discrete setting when errors suffer from long-memory within profiles, but are independent between profiles, and the design points are irregular in both directions. We focus on the case when both $h_1$ and $h_2$ have singularities on the interval $[0, 1]$ but both $1/h_1$ and $1/h_2$ are still integrable. We derive minimax lower  bounds for the $L^2$-risk when $f(t , x)$ belongs to an anisotropic Besov ball and  the blurring function $g(t , x)$ is regular smooth. Combining ideas from the works of Pensky~(2013)  and Chesneau~(2007a), we construct an adaptive wavelet hard-thresholding estimator for $f$ that attains asymptotic minimax convergence rates. Furthermore, we show that such estimator attains convergence rates that depend on the LM parameters and deteriorate as LM strengthens. These rates are completely new and depend on a delicate balance between the smoothness and the spatial homogeneity of the unknown function $f$ in both directions, the degree of ill-posed-ness of the convolution kernel, the long-memory parameter along with the degrees of spatial irregularity associated with $h_1$ and $h_2$. Nonetheless, the spatial irregularity affects convergence rates only if $f$ is spatially inhomogeneous in either direction.  In addition, with $\alpha=1$ and $\beta_1=\beta_2=0$, our rates match exactly those in Benhaddou et al.~(2013) and  Benhaddou~(2017) with $p=2$ and the calibration $\varepsilon^2=\frac{\sigma^2}{MN}$, in their treatment of wavelet regression problem with spatially regular design points and $i.i.d$ Gaussian errors.   Furthermore, if we fix the time variable $t$, our rates are identical to those found in  Antoniadis et al.~(2014), in their treatment of functional deconvolution with spatially irregular design points and $i.i.d$ Gaussian errors as long as data loss is moderate. Finally, if we fix the space variable $x$, our rates match those in Pensky.~(2013), in her treatment of the general spatially inhomogeneous linear inverse problems with singularities and $i.i.d$ Gaussian errors, provided that data loss is moderate.

The rest of the paper is organized as follows. Section 2 introduces some notation as well as the estimation algorithm. Section 3 describes the derivation of the lower bounds for the $L^2$-risk of estimators of $f$ as well as the upper bounds and establishes the asymptotic optimality of the estimator. Finally, Section 4 contains the proofs of the theoretical results.
 \section{Estimation Algorithm.}
 Denote the complex conjugate of $a$ by $\bar{a}$. Consider a finitely supported  periodized wavelet bases (e.g., Daubechies-type), $\psi_{j_1, k_1}(t)$, and $\eta_{j_2, k_2}(x)$. Let $\psi$ be $s_{1o}$-regular  and $\eta$ be $s_{2o}$-regular. Let $m_{10}$ and $m_{20}$ be the lowest resolution levels of $\psi_{j_1, k_1}(t)$ and $\eta_{j_2, k_2}(x)$, respectively and denote the scaling functions for the two bases by  $\psi_{m_{10}-1, k_1}(t)$ and $\eta_{m_{20}-1, k_2}(x)$, respectively. Define
\be \label{eq:setinf}
\Om = \left\{ \bom= (j_1,k_1; j_2,k_2): m_{i0}-1 \leq j_i \leq \infty,\    k_i = 0, \cdots, 2^{j_i}-1; i=1,2 \right\}.
\ee 
 Following Pensky~(2013), find functions $U_{\omega}(t, x)$ such that 
\be \label{Qstaru}
\int^1_0 g(t-s, x)U_{\omega}(s, x)ds=\psi_{j_1k_1}(t)\eta_{j_2k_2}(x).
\ee
Then, applying the Fourier transform in the direction of $t$ to equation \fr{Qstaru} and rearranging yields
\be \label{u1}
U_{\omega}(m, x)= \frac{\overline{\psi_{j_1k_1}(m)}\eta_{j_2k_2}(x)}{g(m,x)},
\ee
where, $g(m, x)$, $\psi_{j_1, k_1}(m)$ are the Fourier coefficients of $g(t, x)$ and $\psi_{j_1k_1}(t)$, respectively. Therefore, 
\be \label{u2}
U_{\omega}(t, x)=\sum_{m\in Z} \frac{\overline{\psi_{j_1k_1}(m)}\eta_{j_2k_2}(x)}{g(m, x)}e^{i2\pi mt}.
\ee
Then, an unbiased estimator analogue to equation (7.6) in Pensky~(2013) is given in its discrete version by 
\be\label{tibeom}
\tilde{\beta}_{\omega}= \frac{1}{MN}\sum^M_{l=1}\sum^N_{i=1}U_{\omega}(t_i, x_l)\frac{{{Y}}(t_i, x_l)}{h_1(t_i)h_2(x_l)}.
\ee
 \noindent
{\bf Assumption A.2.  }\label{A2} The kernel $g(t, x)$ is $\nu-2$ times continuously differentiable in the direction of $t$, and $r>\nu \geq 1$ time differentiable outside of the neighborhood of discontinuity of $g^{(\nu -1)}$ with $g^{(\nu)}$ and $g^{(r)}$ uniformly bounded. In addition, the kernel $g$ is such that the functions $U_{\omega}$ have bounded support in the direction of $t$ and $x$ of lengths proportional to $2^{-j_1}$ and centered at $2^{-j_1}k_1$ in the direction of $t$, and to $2^{-j_2}$ and centered at $2^{-j_2}k_2$ in the direction of $x$. Furthermore, for a fixed $t$, $g$ is uniformly bounded in the direction of $x$, and the functional Fourier coefficients $g(m, x)$, for some positive constants $\nu$, $K_1$ and $K_2$, independent of $m$ and $x$, are such that 
\be  \label{blur}
 K_1\left(|m|+1\right)^{-2\nu} \leq |g(m, x)|^2 \leq K_2 \left(|m|+1\right)^{-2\nu}.
\ee
The parameter $\nu$  represents the degree of ill-posedness in the direction of $t$. \\
\noindent
{\bf Assumption A.3. }\label{A3} The functions $h_1(t)$  and $h_2(x)$ are continuous on the interval $[0, 1]$, and satisfy $h_1(t_o)=0$ and $h_2(x_o)=0$, $t_o, x_o\in(0, 1)$. In addition, there exist some absolute positive constants $c_{h_{i1}}$ and $c_{h_{i2}}$, with $c_{h_{i1}} < c_{h_{i2}}$, such that for any $x$, with $x$, $x+x_o\in[0, 1]$ and $x_o \in (0, 1)$, one has
\be \label{hprop}
c_{h_{i1}}|x_i|^{\beta_i}< h_i(x_i+x_{io})< c_{h{i2}}|x_i|^{\beta_i},\ with\ 0<\beta_i < 1,\ i=1, 2,
\ee
where
\be
x_i=\left\{ \begin{array}{ll} 
t , & \mbox{if}\ \  i=1,\\
x   , & \mbox{if}\ \  i=2,
   \end{array} \right.     ,
\ \ \ x_{io}=\left\{ \begin{array}{ll} 
t_o , & \mbox{if}\ \  i=1,\\
x_o   , & \mbox{if}\ \  i=2.
   \end{array} \right.  
   \ee
Assumption {\bf A.3} corresponds to the situation of moderate data losses in the direction of  $t$ and $x$. Hence, by Assumptions {\bf A.2} and {\bf A.3}, in particular, since $0<\beta_1< 1$ the function $f(\cdot, \cdot)$ can be expanded into a wavelet series as 
\be
f(t, x)= \sum_{\bom \in \Omega}\beta_{\bom}\psi_{j_1, k_1}(t)\eta_{j_2, k_2}(x),
\ee
Choose $m_{i0}$ according to formula $(4.2)$ in Pensky~(2013) and define the sets
\be \label{eq:setOmJ}
\Om(J_1, J_2) = \left\{ \bom= (j_1,k_1; j_2,k_2): m_{i0}-1 \leq j_i \leq J_i-1,\    k_i = 0, \cdots, 2^{j_i-1}; i=1,2  \right\}.
\ee 
Then, allow the hard thresholding estimator for $f(t, x)$
\be \label{ef-hat}
\widehat{f}_{MN}(t, x)= \sum_{\bom \in \Omega(J_1, J_2)}\tilde{\beta}_{\bom} \II \left(|\tilde{\beta}_{\bom}| > \lambda(j_1, j_2) \right)\psi_{j_1, k_1}(t)\eta_{j_2, k_2}(x).
\ee
It remains to determine the choices of  $J_1$, $J_2$ and $\lambda(j_1, j_2)$ in \fr{ef-hat}.
\begin{lemma} \label{lem:U}
Let $U_{\omega}$ be defined in \fr{u2}, and let $h_1$ and $h_2$ satisfy \fr{hprop}. Then, under condition \fr{blur} and provided that $s_{10} > \max\{\nu, r\}$, one has
\beqn 
\int^1_0\int^1_0\frac{U^2_{\omega}(t, x)}{h_1(t)h_2(x)}dtdx &\asymp& 2^{(2\nu+\beta_1) j_1 + \beta_2 j_2}\Pi^2_{i=1}|k_i-k_{i0}|^{-\beta_i},\label{un}\\
\int^1_0\int^1_0\frac{U^4_{\omega}(t, x)}{h^3_1(t)h^3_2(x)}dtdx &\asymp& \frac{2^{j_1(4\nu +3\beta_1) + j_2(3\beta_2+1)}}{\Pi^2_{i=1}|k_i-k_{i0}|^{3\beta_i}}, \label{ufour}
\eeqn
where $k_{01}=t_o2^{j_1}$ and $k_{02}=x_o2^{j_2}$, $j_i=m_{i0}-1, \cdots, {J_i}-1$, $i=1, 2$.
\end{lemma}
\begin{lemma} \label{lem:Var}
Let $\tilde{ \beta}_{\bom}$ be defined in \fr{tibeom} and let $s_{10} > \max\{\nu, r\}$. Then, under the conditions \fr{assum1},  \fr{blur} and \fr{hprop}, as $M, N \rightarrow \infty$, simultaneously, one has
\beqn 
\EE \left|\tilde{ \beta}_{\bom}-\beta_{\bom}\right|^{2} &\asymp&  \frac{\sigma^2}{MN^{\alpha}}2^{(2\nu+\beta_1) j_1 + \beta_2 j_2}\Pi^2_{i=1}|k_i-k_{i0}|^{-\beta_i},\label{var}\\
\EE|\tilde{\beta}_{\bom}-\beta_{\bom}|^4 &\asymp&  \frac{\sigma^4}{M^3N^2}\frac{2^{j_1(4\nu +3\beta_1) + j_2(3\beta_2+1)}}{\Pi^2_{i=1}|k_i-k_{i0}|^{3\beta_i}} +\frac{\sigma^4}{M^2N^{2\alpha}}\frac{2^{2(2\nu +\beta_1) j_1 + 2\beta_2 j_2}}{\Pi^2_{i=1}|k_i-k_{i0}|^{2\beta_i}}. \label{fourmomb}
\eeqn
\end{lemma}
According to Lemma \ref{lem:Var}, choose the thresholds $\lambda(j_1, j_2)$  as  
\be \label{thresh}
\lambda^2_G(j_1,j_2)=\gamma^2 \frac{\sigma^2 2^{(2\nu+\beta_1) j_1 +\beta_2j_2}\ln(MN^{\alpha})}{\Pi^2_{i=1}|k_i-k_{i0}|^{\beta_i}{MN^{\alpha}}}, 
\ee
where $\gamma$ is some positive constant independent of $M$ and $N$, when the errors are Gaussian and

\be \label{thresSG}
\lambda^2_{sG}(j_1,j_2)= \frac{\sigma^2 2^{(2\nu +\beta_1) j_1 + \beta_2j_2}\left[1+\mu^2\ln(MN^{\alpha})\right]}{\Pi^2_{i=1}|k_i-k_{i0}|^{\beta_i}{MN^{\alpha}}},
\ee
when the errors are sub-Gaussian. In addition, the highest resolution levels $J_1$ and $J_2$ should be chosen such that

\be \label{J}
2^{J_1}=\left[\frac{A^2MN^{\alpha}}{\sigma^2}\right]^{\frac{1}{2\nu+1}},\  2^{J_2}=\left[\frac{A^2MN^{\alpha}}{\sigma^2}\right].
\ee 
 \section{Convergence rates and asymptotic optimality.}
Denote 
\beqn  \label{eq10}
 s^{*}_i&=&s_i+1/2 - 1/p,\\
 s''_i&=& \frac{1}{1-\beta_i}\left(\frac{1}{p}-\frac{1}{2}\right),\ \ \ i=1, 2.\\
  s'_i&=&s_i + 1/2 -1/p', \ \ \ p'=\min\{2, p\}.
\eeqn
\noindent
{\bf Assumption A.4.} \label{A4} The function $f(t, x)$ belongs to an anisotropic Besov space. In particular, if $s_o \geq s_2$, its wavelet coefficients $ \beta_{\omega}$ satisfy
\be  \label{eq11}
 B^{s_1, s_2}_{p, q}(A)=\left \{ f \in L^2(U): \left( \sum_{j_1, j_2} 2^{(j_1s_1^{*}+j_2s_2^*)q}\left (\sum_{k_1, k_2}| \beta_{j_1, k_1, j_2, k_2}|^{p}\right)^{q/{p}}\right )^{1/q} \leq A\right \}.
\ee
To construct minimax lower bounds for the $L^2$-risk, we define the $L^2$-risk over the set  $\Theta$ as 
\be  \label{eq12}
 R_{MN}(\Theta)=\inf_{\tilde{f}_{MN}} \sup _{f \in \Theta}\EE \| \tilde{f}_{MN}-f\|^2,
\ee
where $\|g\|$ is the $L^2$-norm of a function $g$ and the infimum is taken over all possible estimators. 
\begin{remark}	
Notice that the quantities \fr{thresh},   \fr{thresSG} and \fr{J} are independent of the parameters of the Besov ball $B^{s_1, s_2}_{p, q}(A)$, and therefore estimator \fr{ef-hat} is adaptive with respect to those parameters. 
 \end{remark}
 \begin{theorem} \label{th:lowerbds}
Let $\min\{s_1, s_2\} \geq \max\{\frac{1}{p}, \frac{1}{2} \}$ with $1 \leq p,q \leq \infty$, and $A > 0$. Then, under conditions \fr{assum1}, \fr{blur}, \fr{hprop} and \fr{eq11}, as $M, N \rightarrow \infty$, simultaneously, one has
 \be \label{lowerbds}
R_{MN} (B^{s_1, s_2}_{p, q}(A)) \geq C A^2\left\{ \begin{array}{ll} 
 \left[\frac{\sigma^2}{A^2MN^{\alpha}} \right]^{\frac{2s_2}{2s_2+1}} , & \mbox{if}\ \  {s_1}> {s_2}(2\nu +1) , s^*_1 > s_2(2\nu +\beta_1), {s_2} > s''_2,\\
 \left[ \frac{\sigma^2}{A^{2}MN^{\alpha}}\right]^{\frac{2s_1}{2s_1 +2\nu +1}}   , & \mbox{if}\ \  s''_1{(2\nu +1)}\leq {s_1}\leq {s_2}{(2\nu +1)}, \beta_2s_1<s^*_2(2\nu+1),\\
  \left[\frac{\sigma^2}{A^2MN^{\alpha}} \right]^{\frac{s^*_1}{2s^*_1+2\nu + \beta_1}}, &  \mbox{if}\  s^*_1 < \frac{s^*_2}{\beta_2}(2\nu + \beta_1), s^*_1 < s_2(2\nu+\beta_1), {s_1}< s''_1{(2\nu +1)},\\
   \left[\frac{\sigma^2}{A^2MN^{\alpha}} \right]^{\frac{s^*_2}{2s^*_2+ \beta_2}}, &  \mbox{if}\  {s^*_1}> \frac{s^*_2}{\beta_2}(2\nu +\beta_1), \ \beta_2s_1 > s^*_2(2\nu +1),
    {s_2} \leq s''_2.
   \end{array} \right.
\ee
  \end{theorem}
\begin{lemma} \label{lem:Lar-D} 
Let $\tilde{\beta}_{\bom}$ and $\lambda(j_1, j_2)$ be defined in \fr{tibeom} and \fr{thresh}, respectively. Let conditions \fr{assum1} and \fr{blur} hold. Then, for some positive constant $\gamma$, as $M$, $N\rightarrow \infty$, one has 
\be \label{pro}
\Pr\left(|\tilde{\beta}_{\bom}-\beta_{\bom}| > \frac{1}{2 }\lambda(j_1, j_2)\right)=O\left(\left[{MN^{\alpha}}\right]^{-\tau}\right).
\ee
\end{lemma}
\begin{theorem} \label{th:upperbds}
Let $\widehat{f}(. , .)$ be the wavelet estimator in \fr{ef-hat}, with $\lambda({j_1, j_2})$ given by \fr{thresh} or \fr{thresSG} and, $J_1$ and $J_2$ given by \fr{J}. Let $\min\{s_1, s_2\} \geq \max\{\frac{1}{p}, \frac{1}{2} \}$, and let conditions \fr{assum1},  \fr{blur}, \fr{hprop} and \fr{eq11} hold. If $\gamma$ in \fr{thresh} or $\mu$ in \fr{thresSG} is large enough, then, as $M, N \rightarrow \infty$, 
 \be \label{upperbds}
R_{MN}( B^{s_1, s_2}_{p, q}(A)) \leq C A^2\left\{ \begin{array}{ll} 
  \left[\frac{\sigma^2\ln(n)}{A^2MN^{\alpha}} \right]^{\frac{2s_2}{2s_2+1}}\left[\ln(n)\right]^{\xi_1} , & \mbox{if}\  {s_1}> {s_2}(2\nu +1) , \frac{s'_1}{s_2} > (2\nu +\beta_1), {s_2} > s''_2,\\
   \left[\frac{ \sigma^{2}\ln(n)}{A^2MN^{\alpha}} \right]^{\frac{2s_1}{2s_1 +2\nu +1}} , & \mbox{if} \  s''_1{(2\nu +1)}\leq {s_1}\leq {s_2}{(2\nu +1)}, s_1<\frac{s'_2}{\beta_2}(2\nu+1),\\
    \left[\frac{\sigma^2\ln(n)}{A^2MN^{\alpha}} \right]^{\frac{2s'_1}{2s'_1+2\nu + \beta_1}}, &  \mbox{if}\  s'_1 < \frac{s'_2}{\beta_2}(2\nu + \beta_1), \frac{s'_1}{s_2} < (2\nu+\beta_1), {s_1}<{s''_1}{(2\nu +1)},\\
     \left[\frac{ \sigma^{2}\ln(n)}{A^2MN^{\alpha}} \right]^{\frac{2s'_2}{2s'_2+\beta_2}}\left[\ln(n)\right]^{\xi_2}, &  \mbox{if}\  {s'_1}\geq \frac{s'_2}{\beta_2}(2\nu +\beta_1), \ s_1 > \frac{s'_2}{\beta_2}(2\nu +1),
    {s_2} \leq s''_2,
\end{array} \right.
\ee
where $\xi_1$ and $\xi_2$ are defined as 
\be \label{xi}
\xi_1 =\left[\II \left( p<2 \right)\II \left( \beta_1= \beta_2 \right) +\II \left( p\geq2 \right)\right]\left( s_1= s_2(2\nu +1) \right), \ \ \
\xi_2 =\II  \left( \beta_2s'_1= s'_2(2\nu + \beta_1) \right),
\ee
and $n=MN^{\alpha}$ is the effective sample size.
\end{theorem}
\begin{remark}	
 Theorems \ref{th:lowerbds} and \ref{th:upperbds} imply that, estimator \fr{ef-hat} is asymptotically near-optimal within a logarithmic factor of $n=MN^{\alpha}$, over a wide range of Besov balls $ B^{s_1, s_2}_{p, q}(A)$. 
 \end{remark}
\begin{remark}	
 The rates of convergence are expressed in terms of the smoothness and the spatial homogeneity of the unknown function $f$, the degree of ill-posed-ness $\nu$ of the convolution operator, the LM parameter $\alpha$ and the degrees of spatial irregularity $\beta_1$ and $\beta_2$ associated with the density functions $h_1$ and $h_2$, respectively. However, the spatial irregularity has a detrimental effect on the convergence rates only if function $f$ is spatially inhomogeneous.
  \end{remark}
\begin{remark}	
  With $\alpha=1$ and $\beta_1=\beta_2=0$, our rates match exactly those in Benhaddou et al.~(2013) and  Benhaddou~(2017) with $p=2$ and the calibration $\varepsilon^2=\frac{\sigma^2}{MN}$, in their treatment of functional deconvolution with spatially regular design points.
    \end{remark}
\begin{remark}	
  If we fix the time variable $t$, our rates are identical to those found in  Antoniadis et al.~(2014), in their treatment of wavelet regression problem with spatially irregular design points as long as data loss is moderate.
    \end{remark}
   \begin{remark}	
  If we fix the space variable $x$, our rates match those in Pensky.~(2013), in her treatment of the general spatially inhomogeneous linear inverse problems with singularities and $i.i.d$ Gaussian errors, provided that data loss is moderate.
  \end{remark}
\begin{remark}	
 Note that the design densities $h_1$ and $h_2$ are assumed to be known, but in practice this may not be the case. These functions can be estimated from the data and their empirical counterpart may be used in formula  \fr{tibeom}, and the interested reader may refer to Kerkyacharian and Picard~(2004). Therefore, estimating $h_1$ and $h_2$ is beyond the scope of this work and we assume that these quantities are known. 
   \end{remark}
\begin{remark}	
    Note that the parameter of LM, $\alpha$, may not be known in advance but can be estimated (see, e.g., Fischer and Akay~(1996), Taqqu, Teverovsky and Willinger~(1997), Pilgram and Kaplan~(1998), and Vivero and Heath~(2012)). Providing completely data driven estimates for  this parameter is beyond the scope of this work and we assume that $\alpha$ is known. In addition, $g$ and therefore $\nu$, may not be known either. If one is interested in taking the uncertainty about $g$ into account, this can be achieved using the methodology of Benhaddou and Liu~(2019) or Hoffmann and Reiss~(2008). 
  \end{remark}
  \begin{remark}	
 Notice that a more general version of condition {\fr{assum1}} is the situation when the level of long-memory differs from one profile $x_l$ to another. In such case, a weighted version of \fr{tibeom}, where the weights depend on the long-memory parameters $\alpha_1, \alpha_2, \cdots, \alpha_M$ and $N$ should achieve asymptotically near-optimal convergence rates with minimal conditions on the number of profiles $M$. 
\end{remark}
\section{Proofs. }
{\bf Proof of Theorem \ref{th:lowerbds}.} \\
\underline{\bf  The dense-dense case.}
Using the same test functions $f_{\tilde{\omega}}$ and $f_{\omega}$ as in Benhaddou et al.~(2013), it can be shown that the $L^2$-norm of the difference satisfies    
\be  \label{Norm}
\|f_{\tilde{\omega}}-f_{\omega}\|^2_2 \geq \rho^2_{j_{1}j_2}  2^{j_1+j_2}/8.
\ee
In order to apply Lemma A.1 of Bunea et al.~(2007), one needs to verify condition $(ii)$. Denote $Q^{(N)}_{l, \omega}$ and $Q^{(N)}_{l, \tilde{\omega}}$, the vectors with components
\beqn
q_{\omega}(t_i, x_l)&=&g(t_i-., x_l)*f_{\omega}(., x_l), \ \ \ i=1, 2, \cdots, N.\\
q_{\tilde{\omega}}(t_i, x_l)&=&g(t_i-., x_l)*f_{\tilde{\omega}}(., x_l), \ \ \ i=1, 2, \cdots, N.
\eeqn
Then, with the help of  {\bf Proposition 2} in Pensky~(2013) adapted to kernel $g$ being bivariate, the Kullback divergence is 
\beqn
K(P_{f_{\omega}}, P_{f_{\tilde{\omega}}})&=&\frac{1}{2\sigma^2}\sum^M_{l=1}\EE_h\left(Q^{(N)}_{l, \omega}- Q^{(N)}_{l, \tilde{\omega}}\right)^T(\Sigma^l_N)^{-1} \left(Q^{(N)}_{l, \omega}- Q^{(N)}_{l, \tilde{\omega}}\right)\nonumber\\
&\leq& \frac{1}{2\sigma^2}\sum^M_{l=1}\lambda_{\max}\left[\left(\Sigma^l_N\right)^{-1}\right]\EE_h\|Q^{(N)}_{l, \omega}- Q^{(N)}_{l, \tilde{\omega}}\|^2\nonumber\\
&=&  \frac{CMN^{\alpha}}{\sigma^2}\rho^2_{j_1, j_2}\sum^{2^{j_1}-1}_{k_1=0}\sum^{2^{j_2}-1}_{k_2=0}\left|\omega_{k_1k_2}-\tilde{\omega}_{k_1k_2}\right|^2\int^1_0\int^1_0\left(g*\psi_{j_1k_1}(t)\right)^2\eta^2_{j_2k_2}(x)h_1(t)h_2(x)dtdx\nonumber\\
&\leq& \frac{CMN^{\alpha}}{\sigma^2}\rho^2_{j_1, j_2}\sum^{2^{j_1}-1}_{k_1=0}2^{-(2\nu +\beta_1) j_1-\beta_2j_2}|k_1-k_{01}|^{\beta_1}\sum^{2^{j_2}-1}_{k_2=0}|k_2-k_{02}|^{\beta_2} \nonumber\\
&\leq& \frac{CMN^{\alpha}}{\sigma^2}\rho^2_{j_1, j_2}2^{-(2\nu +\beta_1) j_1-\beta_2j_2}2^{j_1(\beta_1+1)}2^{j_2(\beta_2+1)}.
\eeqn
Hence, using argument similar to Benhaddou et al.~(2013), the lower bounds are 
\be  \label{delta1}
\delta^2 =C A^2 \left\{ \begin{array}{ll}
\left[\frac{\sigma^{2}}{A^{2}MN^{\alpha}}\right]^{\frac{2s_1}{2s_1+2\nu +1}}, & \mbox{if}\ \ s_1 \leq s_2 (2\nu +1),\\
\left[\frac{\sigma^{2}}{A^{2}MN^{\alpha}}\right]^{\frac{2s_2}{2s_2 +1}},&  \mbox{if}\ \  s_1 > s_2 (2\nu +1).
\end{array} \right.
\ee
\underline{\bf  The sparse-dense case.}
Using the same test functions $f_{\tilde{\omega}}$ and $f_{\omega}$ as in Benhaddou et al.~(2013), and following the same procedure as in the dense-dense case, the lower bounds are 
\be \label{delta2}
\delta^2 = C A^2\left\{ \begin{array}{ll}
\left[\frac{\sigma^{2}}{A^{2}MN^{\alpha}}\right]^{\frac{2s_2}{2s_2 +1}},&  \mbox{if}\ \  s^*_1 \geq s_2 (2\nu +\beta_1) ,\\
\left[\frac{\sigma^{2}}{A^{2}MN^{\alpha}}\right]^{\frac{2s^*_1}{2s^*_1+2\nu +\beta_1}},&  \mbox{if}\ \  \ s^*_1 < s_2 (2\nu+\beta_1).
\end{array} \right.
\ee
\underline{\bf  The dense-sparse case.} Using test functions $f_{\omega}$ defined as
\be
f_{\omega}(t, x)=\rho_{j_1j_2}\sum^{2^{j_1}-1}_{k_1=0}\omega_{k_1}\psi_{j_1k_1}(t)\eta_{j_2k_2}(x),
\ee
and following the same procedure as the dense-dense case, the lower bounds will have the form
\be \label{delta3}
\delta^2 = C A^2\left\{ \begin{array}{ll}
\left[\frac{\sigma^{2}}{A^{2}MN^{\alpha}}\right]^{\frac{2s^*_2}{2s^*_2 +\beta_2}},&  \mbox{if}\ \  \beta_2s_1 \geq s^*_2( 2\nu +1),\\
\left[\frac{\sigma^{2}}{A^{2}MN^{\alpha}}\right]^{\frac{2s_1}{2s_1+2\nu +1}},&  \mbox{if}\ \  \ \beta_2s_1 < s^*_2 (2\nu +1).
\end{array} \right.
\ee
 To complete the proof, notice that the highest of the lower bounds corresponds to 
\be \label{d}
d=\min \left\{ \frac{2s_1}{2s_1 + 2\nu +1}, \frac{2s_2}{2s_2 +1}, \frac{2s^*_1}{2s^*_1+2\nu +\beta_1}, \frac{2s^*_2}{2s^*_2+\beta_2} \right\}.\ \Box
\ee
{\bf Proof of Lemma \ref{lem:U}.} To prove \fr{un}, notice that by \fr{u2}, {\bf Assumptions A.2} and {\bf A.3}, one has 
\beqn
\int^1_0\int^1_0\frac{U^2_{\omega}(t, x)}{h_1(t)h_2(x)}dtdx &
\asymp& \left[h_1(2^{-j_1}k_1)h_2(2^{-j_2}k_2)\right]^{-1}\int^1_0 \int^1_0 {U_{\omega}^2(t, x)}dtdx\nonumber\\
&\asymp& \left[\left|2^{-j_1}k_1-t_0\right|^{\beta_1}\left|2^{-j_2}k_2-x_0\right|^{\beta_2}\right]^{-1}\int^1_0 \int^1_0 {\left[\sum_{m\in Z} \frac{\overline{\psi_{j_1k_1}(m)}\eta_{j_2k_2}(x)}{g(m, x)}e^{i2\pi mt}  \right]^2}dtdx\nonumber\\
&\asymp& \left|2^{-j_1}k_1-t_0\right|^{-\beta_1}\left|2^{-j_2}k_2-x_0\right|^{-\beta_2}\int^1_0 \int^1_0 {\sum_{m\in Z} \frac{{|\psi_{j_1k_1}(m)|^2}|\eta_{j_2k_2}(x)|^2}{|g(m, x)|^2} }dtdx\nonumber\\
&\asymp& \left|2^{-j_1}k_1-t_0\right|^{-\beta_1}\left|2^{-j_2}k_2-x_0\right|^{-\beta_2} {\sum_{m\in Z}\left(|m|+1\right)^{2\nu}{{|\psi_{j_1k_1}(m)|^2}\int^1_0|\eta_{j_2k_2}(x)|^2} }dx\nonumber\\
&\asymp&  \left|2^{-j_1}k_1-t_0\right|^{-\beta_1}\left|2^{-j_2}k_2-x_0\right|^{-\beta_2} {\sum_{m\in Z}\left(|m|^{2\nu}+1\right){{|\psi_{j_1k_1}(m)|^2}} }. \nonumber\\ \label{35}
&\asymp& 2^{j_1\beta_1+j_2\beta_2}  \left|k_1-k_{10}\right|^{-\beta_1}\left|k_2-k_{20}\right|^{-\beta_2}\left[\|\psi\|^2+ 2^{2\nu j_1}\|\psi^{(\nu)}\|^2\right] . \label{35}
\eeqn
In the last line we used the derivative property of Fourier transform applied to the Fourier coefficients $\psi_{j_1k_1}(m)$, namely, $$\psi_{j_1k_1}(m)=(i2\pi m)^{-\nu}2^{j_1\nu}\int^1_02^{j_1/2}\psi^{(\nu)}(2^{j_1}t-k_1)e^{-2\pi mt}dt.$$
The proof of  \fr{ufour} is very similar so we skip it. $\Box$ \\
{\bf Proof of Lemma \ref{lem:Var}.} Note that 
\be \label{norrv}
\tilde{\beta}_{\bom}-\beta_{\bom}=\frac{\sigma}{MN}\sum^M_{l=1}\sum^N_{i=1}\frac{U_{\omega}(t_i, x_l)}{h_1(t_i)h_2(x_l)}\varepsilon_{i, l}.
\ee
Define the $N$-dimensional vector ${\bf U}_{l}$,  with components $U_{\omega, l}=\frac{U_{\omega}(t_i, x_l)}{h_1(t_i)h_2(x_l)}$. Denote the joint distribution of the pair $(t, x)$ by $h$ and the expectation over that joint by $\EE_h$. Indeed, by conditions \fr{assum1}, \fr{blur} and \fr{hprop}, the variance of \fr{norrv} becomes
\beqn
\EE|\tilde{\beta}_{\bom}-\beta_{\bom}|^2&=&\frac{\sigma^2}{M^2N^2}\sum^M_{l=1}\EE\left[U^T_l\left(\varepsilon^{(l)}(\varepsilon^{(l)})^T\right)U_l\right]
\leq \frac{\sigma^2}{MN}\lambda^2_{\max}\left(\Sigma^{(l)}_N\right)\int^1_0 \int^1_0 \frac{U_{\omega}^2(t, x)}{h_1(t)h_2(x)}dtdx\nonumber\\
&\leq& \frac{c_2\sigma^2}{MN^{\alpha}}\int^1_0 \int^1_0 \frac{U_{\omega}^2(t, x)}{h_1(t)h_2(x)}dtdx.\label{35}
\eeqn
To complete the proof use equation  \fr{un}. To prove \fr{fourmomb}, notice that the fourth moment of \fr{norrv} is of order 
\beqn
\EE|\tilde{\beta}_{\bom}-\beta_{\bom}|^4&=& O\left( \frac{\sigma^4}{M^4N^4}\sum^M_{l=1}\EE_{h}\left[\sum^N_{i=1}\frac{U^2_{\omega}}{h_1(t_i)h_2(x_l)}\sum^N_{i=1}\left(\EE|{\varepsilon}_{i, l}|^2\right)\right]^2+ \left[\EE|\tilde{\beta}_{\bom}-\beta_{\bom}|^2\right]^2\right)\nonumber\\
&=& O\left( \frac{\sigma^4}{M^2N^2}\int^1_0\int^1_0 U^4_{\omega}(t, x)h^{-3}_1(t)h_{2}^{-3}(x)dtdx+ \left[\EE|\tilde{\beta}_{\bom}-\beta_{\bom}|^2\right]^2\right). \label{e4}
\eeqn
Now, since $\sum^N_{i=1}\EE|{\varepsilon}_{i, l}|^2=N$, using equations \fr{un} and \fr{ufour} in the last line completes the proof.  $\Box$\\
{\bf Proof of Lemma \ref {lem:Lar-D}.} Using \fr{35}, in the case when the errors are Gaussian, we apply the Gaussian tail inequality to yield \fr{pro}. When the errors are sub-Gaussian, the proof relies on the following lemma 
\begin{lemma} \label{lem:Hans-W} 
(Hanson-Wright inequality). Let $X=(X_1, X_2, \cdots, X_n)$ be a random vector with independent components such that $\EE(X)=0$ and $\|X_i\|_{\phi}\leq K$. Then, for any matrix $B$ and some absolute constant $c_0 > 0$, one has
\be \label{sgaus}
\Pr\left(|X^TBX-\EE[X^TBX]| > t\right) \leq 2 \exp \left\{-c_0\min\left\{\frac{t^2}{K^4||B||^2_F}, \frac{t}{K^2||B||_{sp}}\right\}\right\}.
\ee
\end{lemma}
 Let $V^{(l)}$, $l=1, 2, \cdots, M$, be the $N$-dimensional vectors with elements
\be \label{vl}
V_{i}(x_l)=\frac{U_{\omega}(t_i, x_l)}{h_1(t_i)h_2(x_l)}.
\ee 
Now form the vectors $E$, $Z$ and $V$ such that
\be
E=\left[\pmb{\varepsilon}^{(1)}_N\ \pmb{\varepsilon}^{(2)}_N\ \cdots \pmb{\varepsilon}^{(M)}_N\right]^T,\ \ \   Z=\left[ \pmb{\eta}^{(1)}_N\ \pmb{\eta}^{(2)}_N\ \cdots\ \pmb{\eta}^{(M)}_N\right]^T,\ \ \ V=\left[V^{(1)}\ V^{(2)}\ \cdots\ V^{(M)}\right]^T,
\ee
and define the block diagonal matrices
\be \label{bloc}
 A= \left[ \begin{array}{cccc}\ \ 
A^{(1)}_N& 0  & \cdots & 0  \\
0&  A^{(2)}_N & \cdots& 0 \\
0&0 & \cdots &A^{(M)}_N \\
\end{array} \right].
\ee
Then,
\be \label{devsqd}
|\widehat{\beta}_{\omega}-\beta_{\omega}|^2= \frac{\sigma^2}{M^2{N^2}}\left[\sum^M_{l=1}\sum^N_{i=1}\frac{U_{\omega}(t_i, x_l)}{h_1(t_i)h_2(x_l)}\varepsilon_{i, l}\right]^2=\frac{\sigma^2}{M^2{N^2}}\left[V^TE\right]^2=\frac{\sigma^2}{M^2{N^2}}\left[V^TAZ\right]^2.
\ee
Take $X=Z$ and $B=A^TVV^TA$. Take expectation of \fr{devsqd} yields
\beqns
\EE|\widehat{\beta}_{\omega}-\beta_{\omega}|^2&=& \frac{\sigma^2}{M^2{N^2}}\EE\left[Tr\left(Z^TA^TVV^TAZ\right)\right]\\
&=& \frac{\sigma^2}{M^2{N^2}}\EE\left[Tr\left(A^TVV^TAZZ^T\right)\right]= \frac{\sigma^2}{M^2{N^2}}\EE\left[Tr\left(VV^TAZZ^TA^T\right)\right]\\
&=& \frac{\sigma^2}{M^2{N^2}}\EE_h\left[Tr\left(VV^T\Sigma\right)\right]\leq \lambda_{\max}\left(\Sigma\right)\EE_{h}Tr\left(VV^T\right)=\lambda_{\max}\left(\Sigma\right)\EE_{h}Tr\left(V^TV\right)\\&\leq& \frac{c_2\sigma^2N^{1-\alpha}}{MN} \int^1_0\int^1_0\frac{|U_{\omega}(t, x)|^2}{h_1(t)h_2(x)}dtdx.
\eeqns
Now we need to evaluate the spectral and Frobenius norms of $B$. Note that $rank(B)=1$ and therefore, a conditioning argument on the joint probability density function $h$, gives 
\beqns
\|B\|^2_{F}=\|B\|^2_{sp}&=&\lambda_{\max}(B^TB)=\rho^2_{\max}(A^TVV^TA)\leq \rho^2_{\max}(AA^T)\rho^2_{\max}(VV^T)\\
&=& \rho^2_{\max}(\Sigma)\rho^2_{\max}(VV^T)= \rho^2_{\max}(\Sigma)\rho^2_{\max}(V^TV).
\eeqns
Consequently, applying {\bf Lemma} \ref{lem:Hans-W} with 
\beqns
t= c_2\mu^2MN^{2-\alpha}\ln(MN)\int^1_0\int^1_0\frac{|U_{\omega}(t, x)|^2}{h_1(t)h_2(x)}dtdx,
\eeqns
yields equation $(28)$ with $\tau=\frac{c_0\mu^2}{K}$.
$\Box$\\
{\bf Proof of Theorem \ref{th:upperbds}.} Recall $d$ in \fr{d} and denote
\be \label{chijj}
\chi_{MN}=\frac{\sigma^2\ln(MN^{\alpha})}{A^2MN^{\alpha}}, \ \ 2^{j_{10}}=(\chi_{MN})^{-\frac{d}{2s'_1}}, \ \ 2^{j_{20}}=(\chi_{MN})^{-\frac{d}{2 s'_2}}. 
\ee
Observe that $\mathbb{E} \| \hat{f}_{MN}- f \|^2
\leq\mathbb{E}_1 +\mathbb{E}_2 +\mathbb{E}_3+\mathbb{E}_4$, where
\beqn
\mathbb{E}_1&=&\sum\limits_{k_1=0}^{2^{m_{10}}-1} \sum\limits_{k_2=0}^{2^{m_{20}}-1}\mathrm{Var}( \tilde{\beta}_{m_{10},k_1,m_{20},k_2}),  \\
\mathbb{E}_2&=&\sum\limits_{\bom \in \Omega(J_1,J_2)} \mathbb{E}\left[|\tilde{\beta}_{\omega}- \beta_{\omega} |^2\mathbb{I}(|\tilde{\beta}_{\omega} |>\lambda^{MN}_{j_1}) \right], \\
\mathbb{E}_3&=&\sum\limits_{\bom \in \Omega(J_1,J_2)} \mathbb{E}\left[| \beta_{\omega}|^2\mathbb{I}(|\tilde{\beta}_{\omega} |<\lambda^{MN}_{j_1}) \right],\\
\mathbb{E}_4&=&\left(\sum\limits_{j_1=m_{10}-1}^{J_1-1} \sum\limits_{j_2=J_2}^\infty+\sum\limits_{j_1=J_1}^\infty \sum\limits_{j_2=m_{20}-1}^{J_2-1} +\sum\limits_{j_1=J_1}^\infty \sum\limits_{j_2=J_2}^\infty \right) \sum\limits_{k_1=0}^{2^{j_1}-1} \sum\limits_{k_2=0}^{2^{j_2}-1}  |\beta_{\omega}|^2.
\eeqn
Combining $\mathbb{E}_1$ and $\mathbb{E}_4$, using \fr{var} in $\mathbb{E}_1$ and \fr{eq11} with \fr{J} in $\mathbb{E}_4$, one has 
\beqn \label{ee14}
\mathbb{E}_1+\mathbb{E}_4
&=&O\left(\frac{\sigma^2}{MN^{\alpha}}+ \left(  \sum\limits_{j_2=J_2}^\infty+\sum\limits_{j_1=J_1}^\infty  \right) A^2 2^{-2(js'_1+j's'_2)} \right)=O\left( \left[A^2\chi_{MN}\right]^d\right).
\eeqn
Notice that $\mathbb{E}_2$ and $\mathbb{E}_3$ can be partitioned as $\mathbb{E}_2\leq R_{21} + R_{22}$ and,  $\mathbb{E}_3 \leq R_{31} + R_{32}$, where
\beqn  
R_{21}&=&  \sum_{\bom \in \Om(J_1, J_2)}\EE \left[\left| \tilde{ \beta}_{\bom}-\beta_{\bom}\right|^2 \II \left(   \left| \tilde{ \beta}_{\bom}-\beta_{\bom}  \right| > \frac{1}{2} \lambda^{MN}_{j_1}\right)  \right], \label{r21}\\
R_{22}&=& \sum_{\bom \in \Om(J_1, J_2)}\EE \left[\left| \tilde{ \beta}_{\bom}-\beta_{\bom}\right|^2 \II \left(  \left|  \beta_{\bom} \right| >  \frac{1}{2}\lambda^{MN}_{j_1} \right)  \right],\label{r22}\\
R_{31}&=& \sum_{\bom \in \Om(J_1, J_2)}\left| { \beta}_{\bom}  \right|^2 \Pr \left( \left| \tilde{ \beta}_{\bom}-\beta_{\bom}  \right| > \frac{1}{2} \lambda^{MN}_{j_1} \right),\label{r31}\\
R_{32}&=& \sum_{\bom \in \Om(J_1, J_2)} \left| { \beta}_{\bom}  \right|^2\II \left(  \left|  \beta_{\bom} \right| <  \frac{3}{2}\lambda^{MN}_{j_1} \right).\label{r32}
\eeqn
Combining \fr{r21} and \fr{r31}, using {\bf Lemmas} \ref{lem:Var} and \ref{lem:Lar-D}, with $\tau\geq 4$ and equation \fr{J}, and applying Cauchy-Schwarz inequality, yields
\be
R_{21} + R_{31} = O  \left( A^2\left[ \frac{\sigma^{2}}{MN^{\alpha}}\right]^{{\tau}/{2}-1}  \right)=  O \left( A^2 \left[\frac{\sigma^2}{MN^{\alpha}}\right] \right). \label{r21r31}
\ee
Now, combining \fr{r22} and \fr{r32}, and using \fr{var} and \fr{thresh}, gives
\be  \label{r22r32}
 \Delta=R_{22} + R_{32}=  O \left( \sum_{\bom \in \Om(J_1, J_2)}\min \left\{ \left| { \beta}_{\bom}  \right|^2,   2^{2j\nu} \left[ \chi_{MN}\right]\right\} \right).
\ee
Finally, $\Delta$ can be decomposed into the following components
\beqn  
 \Delta_1&=&  O \left( \left\{\sum^{J_1-1}_{j_1=j_{10}+1}\sum^{J_2-1}_{j_2=m_{0}}+ \sum^{J_1-1}_{j_1=m_{0}}\sum^{J_2-1}_{j_2=j_{20}+1}\right\}A^p 2^{-j_12s'_1} 2^{-j_22s'_2} \right), \label{del1}\\
 \Delta_2&=&O \left( \sum^{j_{10}}_{j_1=m_{0}} \sum^{j_{20}}_{j_2=m_{0}}A^2  2^{j_1(2\nu +1)+j_2}\left[\frac{\sigma^{2}\ln(MN^{\alpha})}{A^2MN^{\alpha}}\right] \II \left(\Xi \right) \right), \label{del2}\\
 \Delta_3&=& O \left( \sum^{j_{10}}_{j_1=m_{0}} \sum^{j_{20}}_{j_2=m_{0}}A^{2\gamma}2^{-\gamma(j_1s'_1+j_2s'_2)} \left[ A^2  [\chi_{MN}] 2^{j_1(2\nu_ +\beta_1)+j_2\beta_2} \right] ^{1-\gamma}\II \left(\Xi^c \right)\right),
 \label{del3}
\eeqn
where $\Xi=\left\{ j_1, j_2: 2^{j_1(2\nu +1)+j_2} \leq   \left[ \chi_{\varepsilon}\right]^{d-1} \right \}$.\\
 \underline{Case 1: ${s_1} \geq {s_2}{(2\nu+1)}$,  ${s'_1} \geq {s'_2}{(2\nu+\beta_1)}$ and  ${(1-\beta_2)s_2} > (\frac{1}{p'}-\frac{1}{p}) $.} In this case, $d=\frac{2s_2}{2s_2+ 1}  $, and
\beqn
			\Delta_3
			&=&O \left(   A^2 [\chi_{MN}]^{1-\gamma}\sum\limits_{j_1=m_{0}}^{j_{10}} 2^{-j_1[\gamma2s'_1-(2\nu + \beta_1)(1-\gamma)]}  \sum\limits_{j_2=m_{0}}^{j_{20}} 2^{-j_2[\gamma2s'_2-\beta_2(1-\gamma) ]} \II\left(  \Xi^c \right)  \right) \nonumber\\
			&=&O \left(  A^2    [\chi_{MN}]^{\frac{2 s_2}{2s_2+1} } \sum\limits_{j_1=m_{0}}^{j_{10}} 2^{  {-j_1}[\gamma2s'_1- (2\nu + \beta_1)(1-\gamma)-(2\nu+1)[\gamma2s'_2-\beta_2(1-\gamma) ]] } \right)\nonumber\\
			&=&O \left(  A^2    [\chi_{MN}]^{\frac{2 s_2}{2s_2+1} } \sum\limits_{j_1=m_{0}}^{j_{10}} 2^{  \frac{-j_1}{d_1}[2(1-\beta_2)(s_1-s_2(2\nu+1))+2(\frac{1}{p'}-\frac{1}{p})(\beta_2-\beta_1)] } \right)\nonumber\\
			&=&O\left( A^2 [\chi_{MN}]^{\frac{2 s_2}{2s_2+1} }    \left[\ln(MN^{\alpha})\right]^{\xi_1} \right),
\eeqn
where $\xi_1$ appears in \fr{xi}, $\gamma=\frac{1-\beta_2}{1-\beta_2+2(\frac{1}{p'}-\frac{1}{p})}$ and $d_1=1-\beta_2 +2(\frac{1}{p'}-\frac{1}{p})$.\\
\underline{Case 2:  $s''_1{(2\nu +1)}< {s_1}< {s_2}{(2\nu +1)}$ and $s_1<\frac{s'_2}{\beta_2}(2\nu+1)$}. In this case, $d=\frac{2 s_1}{2s_1+2\nu+1} $, and
\beqn
	\Delta_3&=&O \left(   A^2 [\chi_{MN}]^{1-\gamma}\sum\limits_{j_1=m_{0}}^{j_{10}} 2^{-j_1[\gamma2s'_1-(2\nu + \beta_1)(1-\gamma)]}  \sum\limits_{j_2=m_{0}}^{j_{20}} 2^{-j_2[\gamma2s'_2-\beta_2(1-\gamma) ]} \II\left(  \Xi^c \right)  \right)\nonumber\\
			         &=&O \left(   A^2    [\chi_{MN}]^{\frac{2 s_1}{2s_1+2\nu+1}} \sum\limits_{j_2=m_{0}}^{j_{20}}2^{-\frac{j_2}{d_2}[2(1-\beta_1)(s_2(2\nu+1)-s_1)+2(\frac{1}{p'}-\frac{1}{p})(2\nu+1)(\beta_1-\beta_2)]   } \right)\nonumber\\
		         &=&O\left(A^2    [\chi_{MN}]^{\frac{2 s_1}{2s_1+2\nu+1}} \right),
\eeqn
where $\gamma=\frac{1-\beta_1}{1-\beta_1+2(\frac{1}{p'}-\frac{1}{p})}$ and $d_2=\left[1-\beta_1 +2(\frac{1}{p'}-\frac{1}{p})\right](2\nu +1)$.\\
\underline{Case 3:  $s'_1 < \frac{s'_2}{\beta_2}(2\nu + \beta_1)$, $\frac{s'_1}{s_2} < (2\nu+\beta_1)$ and ${s_1}<{s''_1}{(2\nu +1)}$}. In this case, $d= \frac{2s'_1}{2s'_1+2\nu +\beta_1}$,  and
\beqn 
			\Delta_3&=&O \left(  A^2 [\chi_{MN}]^{1-\gamma}\sum\limits_{j_1=m_{0}}^{j_{10}} 2^{-j_1[\gamma2s'_1-(2\nu + \beta_1)(1-\gamma)]}  \sum\limits_{j_2=m_{0}}^{j_{20}} 2^{-j_2[\gamma2s'_2-\beta_2(1-\gamma) ]} \II\left(  \Xi^c \right)\right)\nonumber\\
			  &=&O \left(   A^2    [\chi_{MN}]^{\frac{2 s'_1}{2s'_1+2\nu+\beta_1}} \sum\limits_{j_2=m_{0}}^{j_{20}}2^{-\frac{j_2}{d_3}[2(s'_2(2\nu+\beta_1)-\beta_2s'_1)]   } \right)\nonumber\\
			&=&O \left(  A^2    [\chi_{MN}]^{\frac{2 s'_1}{2s'_1+2\nu+\beta_1}} \right), \label{del3b}
			\eeqn
where $\gamma=\frac{2\nu+\beta_1}{2s'_1+2\nu +\beta_1}$ and $d_3=2s'_1+2\nu +\beta_1$.\\
\underline{Case 4:  ${\beta_2s'_1}\geq {s'_2}(2\nu +\beta_1)$, $s_1 > \frac{s'_2}{\beta_2}(2\nu +1)$ and ${s_2} \leq s''_2$}. In this case, $d= \frac{2s'_2}{2s'_2 +\beta_2}$,  and
\beqn 
			\Delta_3&=&O \left(  A^2 [\chi_{MN}]^{1-\gamma}\sum\limits_{j_1=m_{0}}^{j_{10}} 2^{-j_1[\gamma2s'_1-(2\nu + \beta_1)(1-\gamma)]}  \sum\limits_{j_2=m_{0}}^{j_{20}} 2^{-j_2[\gamma2s'_2-\beta_2(1-\gamma) ]} \II\left(  \Xi^c \right)\right)\nonumber\\
			  &=&O \left(   A^2    [\chi_{MN}]^{\frac{2 s'_2}{2s'_2+\beta_2}} \sum\limits_{j_1=m_{0}}^{j_{10}}2^{-\frac{j_1}{d_4}[2(\beta_2s'_1-s'_2(2\nu+\beta_1))]   } \right)\nonumber\\
			&=&O \left(  A^2    [\chi_{MN}]^{\frac{2 s'_2}{2s'_2+\beta_2}} \left[\ln(MN^{\alpha})\right]^{\xi_2} \right), \label{del3b}
			\eeqn
where $\xi_2$ appears in \fr{xi}, $\gamma=\frac{\beta_2}{2s'_2+\beta_2}$ and $d_4=2s'_2+\beta_2$.\\
		Combining the results from \fr{ee14} to \fr{del3b} completes the proof of \fr{upperbds}.
$\Box$


\begin{thebibliography}{99}

\bibitem{abr}
Abramovich, F. and Silverman, B.W. (1998). Wavelet decomposition
approaches to statistical inverse problems. {\it Biometrika}. {\bf
85}, 115--129.

 \bibitem{aenpenpic}
Antoniadis, A., Pensky, M., Sapatinas, T. (2014),
'Nonparametric regression estimation based on spatially inhomogeneous data: minimax global convergence rates and adaptivity',
 {\it ESAIM: PS}, {\bf 18},  1-41.
 
 \bibitem{ben1} Benhaddou, R., Liu, Q. (2019). Anisotropic functional deconvolution with long-memory noise: the case of a multi-parameter fractional Wiener sheet. {\it Journal of Nonparametric Statistics.} {\bf 31(3)}, 567-595.
 
 \bibitem{ben2} Benhaddou, R., Liu, Q. (2019). Minimax adaptive wavelet estimator for the anisotropic functional deconvolution model with unknown kernel. {\it Communications in Statistics- Theory and Methods.} https://doi.org/10.1080/03610926.2019.1617880.
   
 \bibitem{ben2} Benhaddou, R. (2018). Laplace deconvolution with dependent errors: a minimax study. {\it Journal of Nonparametric Statistics.} {\bf 30(4)}, 1032-1048.
 
 \bibitem{benh3}    Benhaddou, R. (2017). On minimax convergence rates under $L^p$-risk for the anisotropic functional deconvolution model. {\it Statistics and Probability Letters}. {\bf 130}, 120-125.
 
  \bibitem{ben3} Benhaddou, R. (2016). Deconvolution model with fractional Gaussian noise: a minimax study. {\it Statistics and Probability Letters}. {\bf 117}, 201-208.
  
\bibitem{ben4} Benhaddou, R., Kulik, R., Pensky, M., Sapatinas, T. (2014). Multichannel deconvolution with long-range dependence: a minimax study. {\it Journal of Statistical Planning and Inference}. {\bf 148}, 1-19.
 
\bibitem{benpenpic}
Benhaddou, R., Pensky, M., Picard, D. (2013).
Anisotropic denoising in functional deconvolution model with dimension-free convergence rates.
 {\it Electron. J.  Stat.} {\bf 7},  1686-1715.
 
  \bibitem{AyTaq}
Cai, T.T. \& Brown, L.D. (1998),
'Wavelet shrinkage for nonequispaced samples',
 {\it The Annals of Statistics,} {\bf 26(5)},  1783-1799.

\bibitem{btw07}
Bunea, F., Tsybakov, A. and Wegkamp, M.H. (2007).
Aggregation for Gaussian regression.
{\it Ann. Statist.} {\bf 35}, 1674--1697.

\bibitem{ches1}
Chesneau, C. (2007a). Regression in random design: a minimax study. {\em Statistics and Probability Letters}. {\bf 77}, 40--53.

\bibitem{ches1}
Chesneau, C. (2007b). Wavelet block thresholding for samples with random design: a minimax approach under $L_p$-risk. {\em Electronic Journal of Statistics}. {\bf 1}, 331--346.

\bibitem{don}
Donoho, D.L. (1995). Nonlinear solution of linear inverse problems
by wavelet-vaguelette decomposition. {\em Applied and Computational
Harmonic Analysis}. {\bf 2}, 101--126.

\bibitem{donrai}
Donoho, D.L. \& Raimondo, M. (2004), 'Translation invariant
deconvolution in a periodic setting', {\it International Journal of
Wavelets, Multiresolution and Information Processing}, {\bf 14},
415--432.

\bibitem{fisc} Fischer, R., Akay, M. (1996), 'A comparison of analytical methods for the study of fractional Brownian motion', {\it Annals of Biomedical Engineering}, {\bf 24 (4)}, 537-543.

\bibitem{hoffr} Hoffmann, M., Reiss, M. (2008), 'Nonlinear Estimation for inverse problems with error in the operator', {\it Annals of Statistics}, {\bf 36}, 310-336.

\bibitem{johnstone}
Johnstone, I.M., Kerkyacharian, G., Picard, D. and Raimondo, M.
(2004). Wavelet deconvolution in a periodic setting. {\it Journal of
the Royal Statistical Society}, Series B. {\bf 66}, 547--573. (with
discussion, 627--657).

\bibitem{kul4} Kerkyacharian, G., Picard, D. (2004). Regression in random design and warped wavelets. {\it Bernoulli}. {\bf 10}, 1053-1105.

\bibitem{kul4} Kulik, R., Sapatinas, T. Wishart, J. R. (2015). Multichannel deconvolution with long-range dependence: Upper bounds on the $L_p$-risk. {\it Applied and Computational Harmonic Analysis}. {\bf 38}, 357-384.

  \bibitem{Paintrtaq}
Painter, S., Paterson, L. \& Beresford, G. (1995),
'On the distribution of the reflection coefficients and seismic amplitudes',
 {\it Geophysics,} {\bf 60(4)},  1187-1194.
 
\bibitem{penpic}
 Pensky, M. (2013).
Spatially inhomogeneous linear inverse problems with possible singularities.
 {\it Ann. Statist.} {\bf 41(5)},  1686-1715.
 
 \bibitem{pensap1}
Pensky, M., Sapatinas, T. (2009),
'Functional deconvolution in a periodic setting: uniform case',
{\it Annals of Statistics},  {\bf 37}, 73--104. 

\bibitem{pensap2}
Pensky, M., Sapatinas, T.  (2010), 
'On convergence rates equivalency and sampling strategies in a functional deconvolution model',
{\it Annals of Statistics}.  {\bf 38}, 1793--1844.

\bibitem{pensap3}
Pensky, M., Sapatinas, T.  (2011),
'Multichannel boxcar deconvolution with growing number of channels',  
{\it Electronic Journal of Statistics},  {\bf 5},  53-82. 

\bibitem{pil} Pilgram, B., Kaplan, D.T. (1998), 'A comparison of estimators for 1/f noise', {\it Physica D}, {\bf 114}, 108-122.
 
 \bibitem{don}
Robinson, E.A. (1999). Seismic Inversion and Deconvolution: Part B: Dual Sensor Technology. {\em Elsevier, Oxford}. 

\bibitem{Taq1} Taqqu, M. S. (1975), 'Weak convergence to fractional Brownian motion and the Rosenblatt process', {\it Probability Theory and Related Fields}, {\bf 31 (4)}, 287-302.

\bibitem{Taq2} Taqqu, M. S., Teverovsky, V., Willinger, W. (1997), 'Estimators for long-range dependence: an empirical study', {\it Fractals}, {\bf 3 (4)}, 785-798.

\bibitem{viv}
Vivero, O., Heath, P. W. (2012),
'A regularized  estimator for long-range dependent processes', 
{\it Automatica}, {\bf 48 (2)},  287-296.

\bibitem{walter}
Walter, G., Shen, X. (1999). Deconvolution using Meyer wavelets.
{\it Journal of Integral Equations and Applications}. {\bf 11},
515--534.

\bibitem{vershyn}
Vershynin, T. (2011),
'Introduction to the non-asymptotic analysis of random matrices'. arXiv:1011.3027 

\bibitem{wan1} Wang, Y. (1997). Minimax estimation via wavelets for indirect long-memory data . {\it Journal of Statistical Planning and Inference}. {\bf 1}, 45-55.

\bibitem{wish} Wishart, J. M. (2013). Wavelet deconvolution in a periodic setting with long-range dependent  errors. {\it Journal of Statistical Planning and Inference}. {\bf 5}, 867-881.

\end{thebibliography}
\end{document}